\documentclass[twoside,11pt,reqno]{amsart}
\usepackage{amsmath,amsthm,amssymb,amstext,amsfonts,amscd}
\usepackage{graphicx}
\usepackage{multirow}
\usepackage{url}
\usepackage{enumerate}
\numberwithin{equation}{section}
\begin{document}
\title[{\bf On Generalization of Bailey's identity involving product of ...}]
{\bf On Generalization of Bailey's identity involving product of generalized hypergeometric series}
\author[{\bf Y. S. Kim and   A. K. Rathie}]{\bf Yong Sup Kim$^{1}$ and  Arjun K. Rathie$^{2}$}
\address{$^{1}$ Yong Sup Kim, Department of Mathematics Education, Wonkwang University, Iksan, Korea }
\email{yspkim@wonkwang.ac.kr}
\address{$^{2}$ Arjun K. Rathie, Department of Mathematics, Central University of Kerala, Kasaragad 671123, Kerala, India }
 \email{akrathie@cukerala.ac.in}

\subjclass[2010]{Primary 33B20, 33C20; Secondary 33B15, 33C15}
\begin{abstract}
The aim of this research paper is to obtain explicit expressions of 
\begin{enumerate}[(i)]
\item $ {}_1F_1 \left[\begin{array}{c} \alpha \\ 2\alpha + i \end{array} ; x \right]. {}_1F_1\left[\begin{array}{c} \beta \\ 2\beta + j \end{array} ; x \right]$\\
    \item ${}_1F_1 \left[\begin{array}{c} \alpha \\ 2\alpha - i \end{array} ; x \right] . {}_1F_1 \left[\begin{array}{c} \beta \\ 2\beta - j \end{array} ; x \right]$ \\and 
    \item ${}_1F_1 \left[\begin{array}{c} \alpha \\ 2\alpha + i \end{array} ; x \right] . {}_1F_1 \left[\begin{array}{c} \beta \\ 2\beta - j \end{array} ; x \right]$
  \end{enumerate}
  in the most general form for any $i,j=0,1,2,\cdots $\\
  
  For $i=j=0$, we recover well known and useful identity due to Bailey. The results are derived with the help of a well known Bailey's formula involving products of generalized hypergeometric series and generalization of Kummer's second transformation formulas available in the literature.  A few interesting   new as well as known special cases have also been given.
  \end{abstract}
  \maketitle
    \section{Introduction} 
    In the theory of hypergeometric series and  generalized hypergeometric series,  classical summation theorems such as those of Gauss, Gauss's second, Kummer and Bailey for the series ${}_2F_1$; Watson, Dixon, Whipple and Saalsh\"utz for the series ${}_3F_2$ and others play a key role. Applications of these classical summation theorems are well known now.
       
    In a very popular, useful and interesting research paper, Bailey [1] had obtained a large number of known as well as new results involving products of generalized hypergeometric series. 
    
    We recall here the following well known and useful transformation formulas due to Kummer[5]
    \begin{equation} \label{Kummer-1} %%%(1.1)
  e^{-x}\, {}_1F_1 \left[\begin{array}{c} \alpha \\ \beta  \end{array}; x\right] = {}_1F_1 \left[\begin{array}{c} \beta-\alpha\, \\
                \beta\,
                \end{array} ; -x  \right]
  \end{equation}
 and
 \begin{equation} \label{Kummer-2} %%%(1.2)
  e^{-\frac{1}{2}x}\, {}_1F_1 \left[\begin{array}{c}\alpha\, \\ 2\alpha\, \end{array} ; x \right]
   = {}_0F_1 \left[\begin{array}{c}       \overline{\hspace{5mm}}\; \\ \alpha +\frac{1}{2}\,
   \end{array};  \frac{x^2}{16} \right].
  \end{equation}
  Bailey[1] derived (1) by employing the classical Gauss's summation theorem[2,4,12,13]:
  \begin{equation}\label{Gauss1}
  {}_2F_1 \left[ \begin{array}{c} a, b \\ c  \end{array}; 1 \right] = \frac{\Gamma(c) \Gamma(c-a-b)}{|\Gamma(c-a) \Gamma(c-b)}
  \end{equation}
  provided $\Re(c-a-b)>0$
  and (2) by using the classical Gauss's second summation theorem[4,8]:
  \begin{equation}\label{Gauss2}
  {}_2F_1 \left[ \begin{array}{c} a, b \\ \frac{1}{2}(a+b+1)  \end{array}; \frac{1}{2} \right]
   =    \frac{\Gamma\left(\frac{1}{2}\right) \Gamma\left(\frac{1}{2}a + \frac{1}{2}b + \frac{1}{2} \right) }{\Gamma\left(\frac{1}{2}a+\frac{1}{2}\right) \Gamma\left(\frac{1}{2}b+\frac{1}{2}\right)}
  \end{equation}   
  In 1998, Rathie and Choi[11] derived the result (1.2) which is equivalent to 
  \begin{equation}
  {}_1F_1 \left[\begin{array}{c}\alpha\, \\ 2\alpha\, \end{array} ; 2x \right]
  = e^x {}_0F_1 \left[\begin{array}{c}       \overline{\hspace{5mm}}\; \\ \alpha +\frac{1}{2}\,
   \end{array};  \frac{x^2}{4} \right]
  \end{equation}
  by utilizing Gauss's second summation theorem (1.4).
  
  On the other hand, from the theory of differential equations, Preece [6] obtained the following very well known and useful identity involving product of the generalized hypergeometric series : 
\begin{equation}
{}_1F_1 \left[\begin{array}{c}\alpha\, \\ 2\alpha\, \end{array} ; x \right]. {}_1F_1 \left[\begin{array}{c}\alpha\, \\ 2\alpha\, \end{array} ; -x \right]
=  {}_1F_2 \left[\begin{array}{c}\alpha\, \\ \alpha+\frac{1}{2}, 2\alpha\, \end{array} ;  \frac{x^2}{4} \right]
\end{equation}

  In 1928, Bailey generalized the Preece's identity (1.6) in the form
    
   \begin{equation}
{}_1F_1 \left[\begin{array}{c}\alpha\, \\ 2\alpha\, \end{array} ; x \right]. {}_1F_1 \left[\begin{array}{c}\beta\, \\ 2\beta\, \end{array} ; -x \right]
= {}_2F_3 \left[\begin{array}{c}\frac{1}{2}(\alpha+\beta),\, \frac{1}{2}(\alpha+\beta+1) \\ \alpha+\frac{1}{2}, \beta+\frac{1}{2}, \alpha+\beta\, \end{array} ;  \frac{x^2}{4} \right]
\end{equation} 
by employing the following classical Watson's summation theorem[2]

\begin{equation}\label{watson}
  {}_3F_2 \left[ \begin{array}{c} a, b,c \\ \frac{1}{2}(a+b+1), 2c  \end{array}; 1 \right]
   =    \frac{\Gamma\left(\frac{1}{2}\right) \Gamma\left(c+\frac{1}{2}\right) \Gamma\left(\frac{1}{2}a + \frac{1}{2}b + \frac{1}{2} \right) \Gamma\left(c-\frac{1}{2}a - \frac{1}{2}b + \frac{1}{2} \right)  }{\Gamma\left(\frac{1}{2}a+\frac{1}{2}\right)\Gamma\left(\frac{1}{2}b+\frac{1}{2}\right) \Gamma\left(c-\frac{1}{2}a+\frac{1}{2}\right)\Gamma\left(c-\frac{1}{2}b+\frac{1}{2}\right)}
  \end{equation} 
  provided $\Re(2c-a-b) >-1.$\\
  
  It is interesting to mention here that the Bailey identity (1.7) reduces to the Preece's identity (1.6) by taking $\beta=\alpha$.\\
  
  Further, if we use Kummer's first transformation (1.1), the Preece's and Bailey's identities can be re-written in the following forms:
   \begin{equation}
\left\{{}_1F_1 \left[\begin{array}{c}\alpha\, \\ 2\alpha\, \end{array} ; x \right]\right\}^2.
=  e^x {}_1F_2 \left[\begin{array}{c}\alpha\, \\ \alpha+\frac{1}{2}, 2\alpha\, \end{array} ;  \frac{x^2}{4} \right]
\end{equation}
and 
\begin{equation}
{}_1F_1 \left[\begin{array}{c}\alpha\, \\ 2\alpha\, \end{array} ; x \right]. {}_1F_1 \left[\begin{array}{c}\beta\, \\ 2\beta\, \end{array} ; x \right]
= e^x \; {}_2F_3 \left[\begin{array}{c}\frac{1}{2}(\alpha+\beta),\, \frac{1}{2}(\alpha+\beta+1) \\ \alpha+\frac{1}{2}, \beta+\frac{1}{2}, \alpha+\beta\, \end{array} ;  \frac{x^2}{4} \right]
\end{equation}

In 1997, Rathie[9] gave a very short proof of Preece's identity (1.9) by using Kummer's second transformation (1.2) and the following Bailey's product formula[1]
\begin{equation}
{}_0F_1 \left[\begin{array}{c}-\, \\ \rho \, \end{array} ; x \right]. {}_0F_1 \left[\begin{array}{c}-\, \\ \sigma \, \end{array} ; x \right]
= {}_2F_3 \left[\begin{array}{c}\frac{1}{2}(\rho+\sigma),\, \frac{1}{2}(\rho+\sigma-1) \\ \rho, \sigma, \rho+\sigma-1\, \end{array} ; 4x \right]
\end{equation}
and obtained the following two results closely related to (1.9)
\begin{align}
{}_1F_1 \left[\begin{array}{c}\alpha\, \\ 2\alpha\, \end{array} ; x \right].{}_1F_1 \left[\begin{array}{c}\alpha\, \\ 2\alpha+1\, \end{array} ; x \right].
& = e^x \left\{{}_1F_2 \left[\begin{array}{c}\alpha\, \\ \alpha+\frac{1}{2}, 2\alpha\, \end{array} ;  \frac{x^2}{4} \right] \right. \nonumber \\
& \qquad  \left. - \frac{x}{2(2\alpha+1)} {}_1F_2 \left[\begin{array}{c}\alpha+1\, \\ \alpha+\frac{3}{2}, 2\alpha+1\, \end{array} ;  \frac{x^2}{4} \right]\right\}
\end{align}
and 
\begin{align}
{}_1F_1 \left[\begin{array}{c}\alpha\, \\ 2\alpha\, \end{array} ; x \right].{}_1F_1 \left[\begin{array}{c}\alpha\, \\ 2\alpha-1\, \end{array} ; x \right].
& = e^x \left\{{}_1F_2 \left[\begin{array}{c}\alpha\, \\ \alpha+\frac{1}{2}, 2\alpha-1\, \end{array} ;  \frac{x^2}{4} \right] \right. \nonumber \\
& \qquad \left. + \frac{x}{2(2\alpha-1)} {}_1F_2 \left[\begin{array}{c}\alpha\, \\ \alpha+\frac{1}{2}, 2\alpha\, \end{array} ;  \frac{x^2}{4} \right]\right\}
\end{align} 
 In 1998, Choi and Rathie established a very short proof of Bailey's result (1.10) by the same technique developed by Rathie and obtained a few results closely related to (1.10).
 
 Very recently, Choi and Rathie [10] established the generalization of Preece's identity (1.9) and obtained explicit expressions of 
\begin{align}
& (i) \;  \left\{{}_1F_1 \left[\begin{array}{c}\alpha\, \\ 2\alpha+i\, \end{array} ; x \right]\right\}^2.\\
& (ii)  \; \left\{{}_1F_1 \left[\begin{array}{c}\alpha\, \\ 2\alpha-i\, \end{array} ; x \right]\right\}^2.\\
& and \nonumber\\
& (iii) \; {}_1F_1 \left[\begin{array}{c}\alpha\, \\ 2\alpha+i\, \end{array} ; x \right].{}_1F_1 \left[\begin{array}{c}\alpha\, \\ 2\alpha-i\, \end{array} ; x \right]
\end{align}
in the most general case for any $i=0,1,2 \cdots $ 

These results have been obtained with the help of the following results recorded in [7].

\begin{equation}
e^{-\frac{1}{2}x} {}_1F_1 \left[\begin{array}{c}\alpha\, \\ 2\alpha+i\, \end{array} ; x \right]
= \sum_{m=0}^i \frac{(-i)_m (2\alpha-1)_m }{(2\alpha+i)_m \left(\alpha-\frac{1}{2}\right)_m \, m!\, 2^{2m}} x^{m} \; {}_0F_1 \left[\begin{array}{c}-\, \\ \alpha+m+\frac{1}{2}\, \end{array} ; \frac{x^2}{16} \right]
\end{equation}
and 
\begin{equation}
e^{-\frac{1}{2}x} {}_1F_1 \left[\begin{array}{c}\alpha\, \\ 2\alpha-i\, \end{array} ; x \right]
= \sum_{m=0}^i \frac{(-1)^m \,(-i)_m (2\alpha-2i-1)_m }{(2\alpha-i)_m \left(\alpha-i-\frac{1}{2}\right)_m \, m!\, 2^{2m}} x^{m} \; {}_0F_1 \left[\begin{array}{c}-\, \\ \alpha+m-i+\frac{1}{2}\, \end{array} ; \frac{x^2}{16} \right]
\end{equation}
for each $i \in N_0$.\\
\textbf{Remark :} For $i=0$, (1.17) or (1.18) reduces to the well-known Kummer's second transformation (1.2). 

The aim of this research paper is to obtain generalization of Bailey's identity (1.10) in the form of three general results.
\begin{enumerate}[(i)]
\item ${}_1F_1 \left[\begin{array}{c}\alpha\, \\ 2\alpha+i\, \end{array} ; x \right]. {}_1F_1 \left[\begin{array}{c}\beta\, \\ 2\beta+j\, \end{array} ; x \right]$
\item ${}_1F_1 \left[\begin{array}{c}\alpha\, \\ 2\alpha+i\, \end{array} ; x \right]. {}_1F_1 \left[\begin{array}{c}\beta\, \\ 2\beta-j\, \end{array} ; x \right]$ and 
\item ${}_1F_1 \left[\begin{array}{c}\alpha\, \\ 2\alpha-i\, \end{array} ; x \right]. {}_1F_1 \left[\begin{array}{c}\beta\, \\ 2\beta-j\, \end{array} ; x \right]$
\end{enumerate}
for each $i,j \in N_0$.

The results are established with the help of the results (1.17), (1.18) and  the Bailey's product   formula (1.11) by  the same technique developed by Rathie[9]. A few known as well as new results have also been given. 
\section{Main Results}
In this section, we present three general formulas asserted in the following theorem. 
\textbf{Theorem 1.} Each of the following formulas hold true for $i,j \in N_0$.
\begin{align}
 & {}_1F_1 \left[\begin{array}{c}\alpha\, \\ 2\alpha+i\, \end{array} ; x \right]. {}_1F_1 \left[\begin{array}{c}\beta\, \\ 2\beta+j\, \end{array} ; x \right] \nonumber \\
 &= e^x \sum_{m=0}^i \sum_{n=0}^j \frac{(-i)_m \, (-j)_n\,(2\alpha-1)_m\,(2\beta-1)_n \, x^{m+n}}{(2\alpha+i)_m\,(2\beta+j)_n\, \left(\alpha-\frac{1}{2}\right)_m\, \left(\beta-\frac{1}{2}\right)_n \,2^{2m+2n} \, m!\, n!} \nonumber\\
 & \quad \times {}_2F_3 \left[\begin{array}{c}\frac{1}{2}(\alpha+\beta+m+n+1),\, \frac{1}{2}(\alpha+\beta+m+n) \\ \alpha+m+\frac{1}{2}, \beta+n+\frac{1}{2}, \alpha+\beta+m+n\, \end{array} ;  \frac{x^2}{4} \right]
\end{align}

\begin{align}
 & {}_1F_1 \left[\begin{array}{c}\alpha\, \\ 2\alpha-i\, \end{array} ; x \right]. {}_1F_1 \left[\begin{array}{c}\beta\, \\ 2\beta-j\, \end{array} ; x \right] \nonumber \\
 &= e^x \sum_{m=0}^i \sum_{n=0}^j \frac{(-1)^{m+n}\,(-i)_m \, (-j)_n\,(2\alpha-2i-1)_m\,(2\beta-2j-1)_n \, x^{m+n}}{(2\alpha-i)_m\,(2\beta-j)_n\, \left(\alpha-i-\frac{1}{2}\right)_m\, \left(\beta-j-\frac{1}{2}\right)_n \,2^{2m+2n} \, m!\, n!} \nonumber\\
 & \quad \times {}_2F_3 \left[\begin{array}{c}\frac{1}{2}(\alpha+\beta+m+n-i-j+1),\, \frac{1}{2}(\alpha+\beta+m+n-i-j) \\ \alpha+m-i+\frac{1}{2}, \beta+n-j+\frac{1}{2}, \alpha+\beta+m+n-i-j\, \end{array} ;  \frac{x^2}{4} \right]
\end{align}
and 
\begin{align}
 & {}_1F_1 \left[\begin{array}{c}\alpha\, \\ 2\alpha+i\, \end{array} ; x \right]. {}_1F_1 \left[\begin{array}{c}\beta\, \\ 2\beta-j\, \end{array} ; x \right] \nonumber \\
 &= e^x \sum_{m=0}^i \sum_{n=0}^j \frac{(-1)^n\,(-i)_m \, (-j)_m\,(2\alpha-1)_m\,(2\beta-2j-1)_n \, x^{m+n}}{(2\alpha+i)_m\,(2\beta-j)_n\, \left(\alpha-\frac{1}{2}\right)_m\, \left(\beta-j-\frac{1}{2}\right)_n \,2^{2m+2n} \, m!\, n!} \nonumber\\
 & \quad \times {}_2F_3 \left[\begin{array}{c}\frac{1}{2}(\alpha+\beta+m+n-j+1),\, \frac{1}{2}(\alpha+\beta+m+n-j) \\ \alpha+m+\frac{1}{2}, \beta+n-j+\frac{1}{2}, \alpha+\beta+m+n-j\, \end{array} ;  \frac{x^2}{4} \right]
\end{align}

\textbf{Proof :} The proof of the results asserted in theorem are quite simple. In order to prove result (2.1), it is sufficient to prove 
\begin{align}
 & e^{-x}\, {}_1F_1 \left[\begin{array}{c}\alpha\, \\ 2\alpha+i\, \end{array} ; x \right]. {}_1F_1 \left[\begin{array}{c}\beta\, \\ 2\beta+j\, \end{array} ; x \right] \nonumber \\
 &=  \sum_{m=0}^i \sum_{n=0}^j \frac{(-i)_m \, (-j)_n\,(2\alpha-1)_m\,(2\beta-1)_n \, x^{m+n}}{(2\alpha+i)_m\,(2\beta+j)_n\, \left(\alpha-\frac{1}{2}\right)_m\, \left(\beta-\frac{1}{2}\right)_n \,2^{2m+2n} \, m!\, n!} \nonumber\\
 & \quad \times {}_1F_3 \left[\begin{array}{c}\frac{1}{2}(\alpha+\beta+m+n+1),\, \frac{1}{2}(\alpha+\beta+m+n) \\ \alpha+m+\frac{1}{2}, \beta+n+\frac{1}{2}, \alpha+\beta+m+n\, \end{array} ;  \frac{x^2}{4} \right]
\end{align}
Now, denoting left hand side of (2.4) by S  and writing in the following form,
\begin{equation}
S = \left\{e^{-\frac{1}{2}x}\, {}_1F_1 \left[\begin{array}{c}\alpha\, \\ 2\alpha+i\, \end{array} ; x \right]\right\}. \left\{ e^{-\frac{1}{2}x}\,{}_1F_1 \left[\begin{array}{c}\beta\, \\ 2\beta+j\, \end{array} ; x \right] \right\}
\end{equation}
Then using the known result (1.17) to each term of (2.5), we get after some simplification
\begin{align}
S & = \sum_{m=0}^i \sum_{n=0}^j \frac{(-i)_m \, (-j)_n\,(2\alpha-1)_m\,(2\beta-1)_n \, x^{m+n}}{(2\alpha+i)_m\,(2\beta+j)_n\, \left(\alpha-\frac{1}{2}\right)_m\, \left(\beta-\frac{1}{2}\right)_n \,2^{2m+2n} \, m!\, n!}\nonumber\\
& \quad \times\, {}_0F_1 \left[\begin{array}{c}-\, \\ \alpha+m+\frac{1}{2}\, \end{array} ; \frac{x^2}{16} \right]. {}_0F_1 \left[\begin{array}{c}-\, \\ \beta+n+\frac{1}{2}\, \end{array} ; \frac{x^2}{16} \right]
\end{align}
Finally, using Bailey's product formula (1.11), it is very easy to see that the S is equal to the right-hand side of (2.1).\\
This completes the proof of (2.1).\\
In exactly the same manner, the result (2.2) and (2.3) can be established. So the details are omitted. 
\section{Special Cases}
 In (2.1), (2.2) and (2.3), if we set $\beta=\alpha$, we get the following interesting results which are also of general nature. These are
 
\begin{align}
 & {}_1F_1 \left[\begin{array}{c}\alpha\, \\ 2\alpha+i\, \end{array} ; x \right]. {}_1F_1 \left[\begin{array}{c}\alpha\, \\ 2\alpha+j\, \end{array} ; x \right] \nonumber \\
 &= e^x \sum_{m=0}^i \sum_{n=0}^j \frac{(-i)_m \, (-j)_n\,(2\alpha-1)_m\,(2\alpha-1)_n \, x^{m+n}}{(2\alpha+i)_m\,(2\alpha+j)_n\, \left(\alpha-\frac{1}{2}\right)_m\, \left(\alpha-\frac{1}{2}\right)_n \,2^{2m+2n} \, m!\, n!} \nonumber\\
 & \quad \times {}_2F_3 \left[\begin{array}{c}\frac{1}{2}(2\alpha+m+n+1),\, \frac{1}{2}(2\alpha+m+n) \\ \alpha+m+\frac{1}{2},\alpha+n+\frac{1}{2}, 2\alpha+m+n\, \end{array} ;  \frac{x^2}{4} \right]
\end{align}

\begin{align}
 & {}_1F_1 \left[\begin{array}{c}\alpha\, \\ 2\alpha-i\, \end{array} ; x \right]. {}_1F_1 \left[\begin{array}{c}\alpha\, \\ 2\alpha-j\, \end{array} ; x \right] \nonumber \\
 &= e^x \sum_{m=0}^i \sum_{n=0}^j \frac{(-1)^{m+n}\,(-i)_m \, (-j)_n\,(2\alpha-2i-1)_m\,(2\alpha-2j-1)_n \, x^{m+n}}{(2\alpha-i)_m\,(2\alpha-j)_n\, \left(\alpha-i-\frac{1}{2}\right)_m\, \left(\alpha-j-\frac{1}{2}\right)_n \,2^{2m+2n} \, m!\, n!} \nonumber\\
 & \quad \times {}_2F_3 \left[\begin{array}{c}\frac{1}{2}(2\alpha+m+n-i-j+1),\, \frac{1}{2}(2\alpha+m+n-i-j) \\ \alpha+m-i+\frac{1}{2}, \alpha+n-j+\frac{1}{2}, 2\alpha+m+n-i-j\, \end{array} ;  \frac{x^2}{4} \right]
\end{align}
and 
\begin{align}
 & {}_1F_1 \left[\begin{array}{c}\alpha\, \\ 2\alpha+i\, \end{array} ; x \right]. {}_1F_1 \left[\begin{array}{c}\alpha\, \\ 2\alpha-j\, \end{array} ; x \right] \nonumber \\
 &= e^x \sum_{m=0}^i \sum_{n=0}^j \frac{(-1)^n\,(-i)_m \, (-j)_m\,(2\alpha-1)_m\,(2\alpha-2j-1)_n \, x^{m+n}}{(2\alpha+i)_m\,(2\alpha-j)_n\, \left(\alpha-\frac{1}{2}\right)_m\, \left(\alpha-j-\frac{1}{2}\right)_n \,2^{2m+2n} \, m!\, n!} \nonumber\\
 & \quad \times {}_2F_3 \left[\begin{array}{c}\frac{1}{2}(2\alpha+m+n-j+1),\, \frac{1}{2}(2\alpha+m+n-j) \\ \alpha+m+\frac{1}{2}, \alpha+n-j+\frac{1}{2}, 2\alpha+m+n-j\, \end{array} ;  \frac{x^2}{4} \right]
\end{align}
Further in (3.1), (3.2) and (3.3), if we set $j=i$, we immediately recover the results recently obtained by Choi and Rathie[3].\\

\textbf{Concluding Remark : } 
In this paper, we have obtained explicit expressions of 
\begin{enumerate}[(i)]
\item $ {}_1F_1 \left[\begin{array}{c} \alpha \\ 2\alpha + i \end{array} ; x \right]. {}_1F_1\left[\begin{array}{c} \beta \\ 2\beta + j \end{array} ; x \right]$\\
    \item ${}_1F_1 \left[\begin{array}{c} \alpha \\ 2\alpha - i \end{array} ; x \right] . {}_1F_1 \left[\begin{array}{c} \beta \\ 2\beta - j \end{array} ; x \right]$\\ and \\
    \item ${}_1F_1 \left[\begin{array}{c} \alpha \\ 2\alpha + i \end{array} ; x \right] . {}_1F_1 \left[\begin{array}{c} \beta \\ 2\beta - j \end{array} ; x \right]$
  \end{enumerate}
  in the most general form for any $i,j \in N $\\
  
  Applications of these results in obtaining a large number of results(reduction and transformation formulas) involving Exton's and Srivastava's series and also for multiple series have also been obtained. The results will be published soon.
%%=======================================================  
  \renewcommand{\refname}{

\end{document}